# PRESERVICE MATHEMATICS TEACHERS UNDERSTANDING OF MODE

| Md Amiruzzaman | Karl W. Kosko | Stefanie R. Amiruzzaman |
|---|---|---|
| Kent State University | Kent State University | Kent State University |
| mamiruzz@kent.edu | kkosko1@kent.edu | sbodenmi@kent.edu |

*This paper explores two preservice mathematics teachers' understanding of mode. Participants' initial understanding and understanding following use of an interactive virtual manipulative is examined. Findings suggest that participants initially operated with less-effective definitions of mode. Preservice teachers' developed understanding following a learning trajectory is discussed.*

## Purpose

Mode is a statistical representation, often in the form of a single value or classification, for the distribution of a data set. Mode is an important concept that is typically introduced to students as early as middle school. Unfortunately, researchers have found that a significant portion of students at various levels (i.e., middle-grade, high school, and undergraduate) have difficulty understanding and explaining the concept mode (Groth & Bergner, 2006; Jacobbe & Carvalho, 2011). Given that teachers' mathematical conceptions can influence those of their students (Ball, Thames & Phelps, 2008; Jacobbe & Carvalho, 2011), the present study focuses on the nature of teachers' understanding of mode. Thus, it is the purpose of the present study to explore preservice mathematics teachers' (PSMTs) conceptions of mode.

## Perspective

Relatively few studies have focused specifically on conceptions of mode of either teachers or students. Those that have generally done so in the context of or in relation to statistical mean. Groth and Bergner (2006) interviewed PMSTs to investigate their understanding of mean, median, and mode. Their research study suggested that most PSMTs (34 out of 45) believed that a data set may have only one mode, as they defined mode as 'most frequent number.' Only a few participants (3 out of 45) believed that a data set may have more than one mode or no mode at all (Groth & Bergner, 2006). Jacobbe and Carvalho (2011) reported that often preservice and inservice mathematics teachers confuse the concept of the mean with the concept of mode. In fact, many mathematic teachers (over 30%) defined mode incorrectly (Jacobbe & Carvalho, 2011). Barr's (1980) observation of college students revealed that 68% of the participants considered the frequency count, and not the classification (in this case, a number) as the mode. Thus, even at the collegiate level where students learn to become teachers, many consider mode in a manner that does not meet the mathematical definition.

While previous studies reported lack of PSMTs' knowledge about the mode, they only reported PSMTs' amount of understanding, but did not provide reasons as why that is the case. Furthermore, prior study has not examined how conceptions of mode develop or evolve. By contrast, the present research study sought to explore PSMTs' understanding of mode and provide descriptions of learning trajectories regarding conceptions of the mode.

## Method

An individual teaching experiment was used to guide this study. In an *individual teaching experiment*, a researcher (a.k.a. teacher-researcher) interviews individual students on one-on-one basis with the intent to model their conceptual understanding of a particular topic (Steffe & Thompson, 2000). Also, the researcher observes and identifies individuals' actions so that the researcher can model learning processes, conceptions, and less-effective understandings. Because of its flexibility, an individual teaching experiment allows a researcher to focus on one student or participant at a time.

The present study consists of nine teaching episodes across a period of nine weeks. Consistent with teaching experiment design guidelines, an observer-researcher was present in all episodes to improve reliability of researcher-constructed models of PSMTs' conceptions. For each episode, tasks were prepared ahead of time. While main tasks were designed to focus on a specific topic or concept (e.g., definition of mode, finding mode of a given data set), probing questions were used during each episode) to clarify participants' answers and press them to explain or explore aspects of the task.

In order to facilitate PSMTs' engagement in tasks related to the mode, a computer-based virtual manipulative in the form of a number line was used. Number lines are widely considered as an appropriate tool for investigating various measures of central tendency (Amiruzzaman & Kosko, 2016; Gravemeijer, 2004). Thus, the Interactive Statistical Number Line (ISNL) was incorporated within specific tasks (see Figures 1 and 2 for screenshots of ISNL). ISNL allowed PSMTs to manipulate the location of data as discrete elements within a dataset. For mode-based tasks, this generally resembled a line plot. Although visually similar to a paper-based line plot, use of ISNL was hypothesized to engage participants in considering how all elements are represented aspects of the data set, given the virtual manipulation of each element in the data set.

## Analysis and Findings

Two Preservice Mathematics Teachers participated in this study: Alex and Bob. Considering the length of this paper only Alex's answers are presented here. Data were collected from a larger study focusing on various aspects of measures of central tendency (mean, mode, & median). However, only data focusing on PSMTs' interaction with mode-based tasks is reported. Considering the length of this paper, we limit our description of the analysis and findings of Alex's actions within episode 4.

Episode 4 was intended to explore Alex and Bob's understanding of the mode. Given the initial task to find the mode of a given data set, Alex sorted the data set as, 1, 3, 3, 4, 5, 6, 9, and then told me that mode is 3 (see Figure 2). Following this, he used the ISNL to develop a mathematical model for the data set and indicated 3 as a mode of the data set. To explain his work, he said that he needed to sort the data set, so that he could see the repeated numbers. In the example, there are two 3s, Alex saw them together once he sorted the numbers. He said, "[after sorting] …now I see two 3s…so 3 is mode." Alex knew that the mode is a repeated number and his algorithmic scheme helped him to find the mode (see Figure 1). We refer to this scheme as *algorithmic* because Alex followed a step-by-step procedure to find the mode.

To confirm the model and his initial understanding, Alex was asked to define the mode. Alex answered, "What occurs the most is mode, in the example, 3 occurred the most, so 3 is mode." So, Alex was asked to find mode from, 1, 2, 3, 4. He responded that all of them are mode as all of them occurred most. Thus, Alex's algorithmic scheme aligns with a definition of mode as "the most."

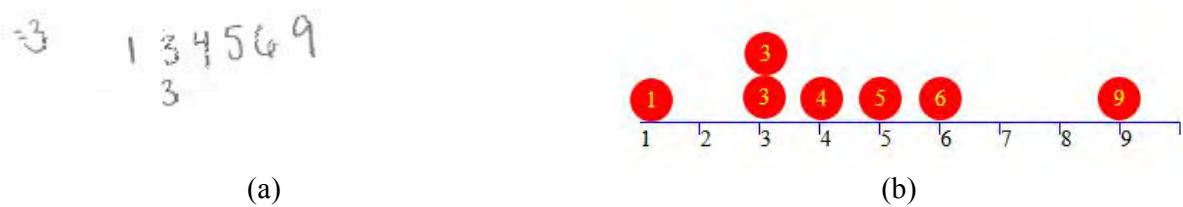

(a)                                                  (b)

**Figure 2.** Alex's approach to find mode from a data set for Task 1 (a) and Task 2 (b).

To press Alex's use of this definition, the teacher-researcher asked him to compare the first data set (1, 3, 3, 4, 5, 6, 9) with the second data set (1, 2, 3, 4). Alex used the ISNL to model both data sets (see Figure 2). While he was comparing models of both data sets, he identified some differences, "I see that … data sets are not same.…in the first data set, there was a repeated number which was the mode… in the second data set there were no repeated numbers…" After analyzing both data sets via the ISNL, he decided to revise his definition of mode (see Figure 2).

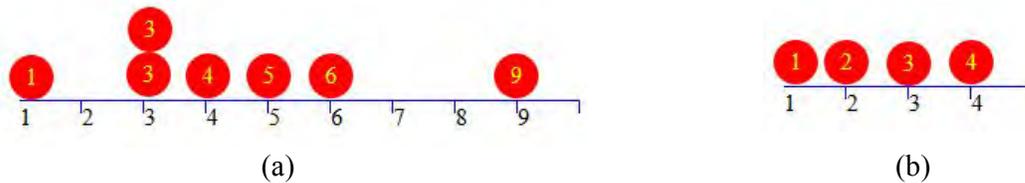

(a)                                                                 (b)

**Figure 3.** Alex's models for two data set. One with repeated number and another without repeated number.

Alex modeled additional datasets on the ISNL to find a mode before revising his definition. He said, "The mode is the value that appears most often in a set of data. A data element will be considered to be a mode if it has a frequency of at least two and there is no other data element that has a higher frequency than that. A data set may have no mode, or one mode, or more than one mode." After revising the definition of mode, Alex went back to find the mode of the data set consists of 1, 2, 3, 4. He analyzed the data set once again and developed a mathematical model using ISNL, and he claimed that the data set 1, 2, 3, 4 has 'no mode'. He further added, "…If there are no data elements that occur most frequently, then that situation indicates that the data set has no mode…"

Alex's initial understanding confirms the finding of Groth and Bergner (2006). Knowing that mode is the most frequent number is not enough. In fact, in some cases 'most frequent number' will help to find a mode of one data set, but not necessarily another (i.e., the first data set versus second). By comparing the two sets on the ISNL, Alex identified that his definition of the mode was incomplete and was not suitable for all scenarios or data sets. We conjecture that by manipulating the individual elements of data along a number line, with the purpose of finding the mode (to represent the set of data). Note that, Alex was confronted with scenarios that were at odds with his prior definitions. The comparison task, in particular, perturbed his acceptance of both his prior definition and operational scheme. Specifically, Alex initially considered his algorithm before considering the nature of the data. By the end of the episode, Alex demonstrated attention to the elements of the data as part of a dataset. Although seemingly a subtle distinction, this shift in how Alex operated on the data to find mode allowed for more flexible considerations, and more correct identifications of mode.

## Conclusions

This paper described one PSMT's (Alex's) points of cognitive dissonance related to the concept of mode, and provides an initial learning trajectory. A key facilitator of Alex's developing conception of mode appeared to be his use of visualization of mathematical models. Such models may be helpful for both PSMTs and K-12 students to develop more sophisticated schemes in operating with data.

It is a common assumption in mathematics education that merely memorizing definitions and working with a few examples is not sufficient to develop a deeper understanding of concepts. Yet mode is often considered a simple concept, and this study provides evidence that many PSMTs (and potential their future students) may not know that they have understood less useful definition of Mathematical Knowledge for Teaching 628mode (or even one considered incorrect by the discipline). Alex considered his definition of mode as insufficient only after he developed mathematical models and visualized both models side-by-side. However, some individuals may have needed additional models that can further press their definition of mode. Specifically, Alex visualized a comparison of multiple datasets, and manipulated elements within those datasets through construction of mathematical models. Such engagement pressed Alex to consider a different way of determining mode, and allowed a transition from an algorithmic scheme to a whole set unit scheme. In the *whole set unit* scheme, Alex considered the whole data set as a unit before finding the mode. A similar transition may be possible with other PSMTs or with middle school students. However, future study is needed to both confirm and extend the findings presented here.

**References**
Amiruzzaman, M., & Kosko, K. W. (2016). Exploring students' understanding of median. Brief research report. In MB Wood, EE Turner, M. Civil, & JA Eli (Eds.), Proceedings of the 38th annual meeting of the North American Chapter of the International Group for the Psychology of Mathematics Education.
Ball, D. L., Thames, M. H., & Phelps, G. (2008). Content knowledge for teaching what makes it special? *Journal of teacher education*, *59*(5), 389-407.
Barr, G. V. (1980). Some student ideas on the median and the mode. *Teaching Statistics, 2*(2), 38-41.
Begg, A., & Edwards, R. (1999). *Teachers' ideas about teaching statistics*. Proceedings of the 1999 Combined Conference of the Australian Association for Research in Education and the New Zealand Association for Research in Education, Melbourne: Australian Association for Research in Education.
Gravemeijer, K. (2004). Local instruction theories as means of support for teachers in reform mathematics education. *Mathematical Thinking and Learning*, *6*(2), 105-128.
Groth, R. E., & Bergner, J. A. (2006). Preservice elementary teachers' conceptual and procedural knowledge of mean, median, and mode. Mathematical Thinking and Learning, 8(1), 37-63.
Jacobbe, T., & Carvalho, C. (2011). Teachers' understanding of averages. In *Teaching Statistics in School Mathematics-Challenges for Teaching and Teacher Education* (pp. 199-209). Springer Netherlands.
Steffe, L. P. & Thompson, P. W. (2000). *Radical Constructivism in Action: Building on the Pioneering Work of Ernst von Glasersfeld*. New York: Falmer.Galindo, E., & Newton, J., (Eds.). (2017). *Proceedings of the 39th annual meeting of the North American Chapter of the International Group for the Psychology of Mathematics Education.* Indianapolis, IN: Hoosier Association of Mathematics Teacher Educators.